\documentclass[11pt]{amsart}
\usepackage[a4paper,margin=28mm]{geometry}

\usepackage{xcolor}
\usepackage{lineno}
\usepackage{amsmath,amssymb,mathrsfs}
\usepackage[shortlabels]{enumitem}
\usepackage{microtype}

\usepackage[
    pagebackref=true,
    colorlinks=false,
    pdfborder={0 0 1},
    linkbordercolor={1 0 0},
    citebordercolor={0 1 0},
    urlbordercolor={0 1 1}
]{hyperref}

\newtheorem{theorem}{Theorem}[section]
\newtheorem{lemma}[theorem]{Lemma}

\newtheorem{corollary}[theorem]{Corollary}
\newtheorem{observation}[theorem]{Observation}

\newtheorem{question}[theorem]{Question}
\theoremstyle{definition}

\theoremstyle{remark}
\newtheorem{remark}[theorem]{Remark}
\numberwithin{equation}{section}
\setlist[enumerate]{label=\textup{(\roman*)},leftmargin=*,itemsep=3pt}
\title[Vanishing on Hartogs triangles]
{Vanishing theorems for $L^2$ Dolbeault cohomology groups of projective Hartogs  triangles}
\author{Jian Chen}
\address{School of Mathematics and Statistics, Central China Normal University,
Wuhan 430079, People's Republic of China}
\email{jian-chen@whu.edu.cn}
\thanks{The author is supported by  NSFC  (Grant No. 12601144).}

\date{\today}

\begin{document}
\begin{abstract}
In this note, using J.-P. Demailly's $L^2$ existence theorem for possibly incomplete K\"ahler metrics, we establish $L^2$ estimates and vanishing theorems for a class of domains in $\mathbb P^n$ without  any boundary regularity assumptions.  In particular, we give an affirmative answer to M.-C. Shaw's question concerning the vanishing of the $L^2$ Dolbeault cohomology groups
of projective Hartogs  triangles.
\end{abstract}
\maketitle

\setcounter{tocdepth}{2}
\tableofcontents

\section{Introduction}\label{sec-introduction}

Throughout this paper, let $\omega_{\rm FS}$ denote the normalized
Fubini--Study form on $\mathbb P^n$ given in standard affine
coordinates $\zeta=(\zeta_1,\ldots,\zeta_n)$ by
\begin{equation}\label{eq-fubini-study}
 \omega_{\rm FS}
 =\sqrt{-1}\,\partial\bar\partial\log(1+|\zeta|^2);
\end{equation}
for differential  forms on a domain $\Omega\subset\mathbb P^n$, all unweighted pointwise norms or
$L^2$ norms use the restricted Fubini--Study metric;
 $\bar\partial^{p,q}$ denotes the
maximal realization of $\bar\partial$ on $(p,q)$-forms:
its domain consists of the square-integrable $(p,q)$-forms
whose $\bar\partial$ (in the sense of distribution) is square integrable; the \emph{$L^2$ Dolbeault cohomology group} is the quotient $H^{p,q}_{(2)}(\Omega)
 =
 \frac{\ker(\bar\partial^{p,q})}
 {\operatorname{im}(\bar\partial^{p,q-1})}.$

For $n\geq 2$, we define the \emph{projective Hartogs triangles} in
$\mathbb P^n$ (see \cite[\S 4]{Sh25} for the case $n=2$) by
\begin{equation}\label{def-intro-tri}
    \begin{aligned}
 \mathbb H_n^+&=
 \{[Z_0:\cdots:Z_n]\in\mathbb P^n:
       |Z_{n-1}|<|Z_n|\},\\
 \mathbb H_n^-&=
 \{[Z_0:\cdots:Z_n]\in\mathbb P^n:
       |Z_n|<|Z_{n-1}|\}.
 \end{aligned}
\end{equation}
Note that  on $\mathbb H^+$,  $Z_n\neq 0$.
Then the map 
$$\mathbb H^+\longrightarrow \mathbb C^{n-1}\times\mathbb D,
\qquad
[Z_0:\cdots:Z_n]\longmapsto
\left(
\frac{Z_0}{Z_n},\ldots,\frac{Z_{n-2}}{Z_n},
\frac{Z_{n-1}}{Z_n}
\right)$$
gives the  biholomorphism. Thus
each Hartogs triangle is biholomorphic to
$\mathbb C^{n-1}\times\mathbb D$
with $\mathbb D$ the unit disk, and is therefore Stein.
In particular, each Hartogs triangle admits a complete
K\"ahler metric (e.g., \cite[Chapter VIII, Theorem 5.2]{Dem12}),
although the restricted Fubini--Study metric is incomplete.

Shaw posed the following question concerning the 
$L^2$ Dolbeault cohomology of projective
Hartogs triangles in $\mathbb P^2$ (\cite{Sh15,Sh23,Sh25}).

\begin{question}[{Shaw's question \cite[Section 6, Problem 5]{Sh25}}]
\label{ques-shaw-intro}
Let $\mathbb H_2^\pm\subset\mathbb P^2$ be the projective Hartogs
triangles defined above, equipped with the restricted
Fubini--Study metric.
Determine whether the $L^2$ Dolbeault cohomology
$H^{1,1}_{(2)}(\mathbb H_2^\pm)$ satisfies the following:
\begin{enumerate}[\rm{(}1\rm{)}]
\item $H^{1,1}_{(2)}(\mathbb H_2^\pm)$ is Hausdorff;
\item $H^{1,1}_{(2)}(\mathbb H_2^\pm)=0$.  \footnote{The vanishing question for  $H^{1,1}_{(2)}(\mathbb H_2^\pm)$. was previously posed in
\cite[p. 113]{Sh15} and repeated in
\cite[Section 4.1, Question (2)]{Sh23}.}
\end{enumerate}
\end{question}

For pseudoconvex domains in $\mathbb P^n$, Shaw
\cite[Theorem 2.9]{Sh25} (see also \cite{BC00,FS22,HI00}  for a proof)  stated an $L^2$ Dolbeault cohomology
vanishing criterion  
\footnote{Shaw further pointed out (\cite[Paragraph after Proposition 2.10]{Sh25}) that the Lipschitz boundary hypothesis
cannot in general be omitted when $p>0$.  Indeed, $H^{2,1}_{(2)}(\mathbb H_2^\pm)$ is even non-Hausdorff (\cite[Theorem 5.1 and Corollary 5.2]{Sh25}).}
 under the additional assumption that
the boundary is Lipschitz.
Chakrabarti--Harrington--Raich
\cite[Theorems 1.4 and 1.5]{CHR26} also established
vanishing criteria for $L^2$ Dolbeault cohomology on
smoothly bounded relatively compact domains in a general complex
manifolds,  under additional assumptions on the complex
Hessian of a smooth weight and the boundary Levi form.   

Note that neither the criterion in \cite[Theorem 2.9]{Sh25} nor those
in \cite[Theorems 1.4 and 1.5]{CHR26} apply directly to the
projective Hartogs triangles $\mathbb H_2^\pm$ to answer
Question \ref{ques-shaw-intro}. Indeed, their boundary is not
Lipschitz at the vertex, while it is Levi-flat at each regular
point and singular at the vertex. Thus the boundary conditions
required in these results are not satisfied.

To address  Question \ref{ques-shaw-intro}, the main tool we
will adopt is Demailly's $L^2$ existence theorem
\cite[Theorem 4.1 and Remark 4.2]{Dem82}
(see Lemma \ref{lem-demailly-existence}), which is
formulated for vector bundle valued $(n,q)$-forms and  requires no boundary regularity. 
To apply Demailly's  theorem, we consider $(p,q)$-forms as
$(K_{\mathbb P^n}^{-1}\otimes
\Omega_{\mathbb P^n}^p)$-valued $(n,q)$-forms.  The main difficulty  is that 
the  Hermitian metric (induced by Fubini--Study metric $\omega_{\rm FS}$) on this
bundle is only Nakano semipositive
(see Lemma \ref{lem-curvature}).
Our strategy is to obtain a uniform positive Nakano
curvature lower bound on $\mathbb H_n^\pm$ by
introducing a suitable weight.

Note that each projective Hartogs triangle $\mathbb H_n^\pm$
admits a bounded smooth plurisubharmonic function whose
Levi form has an eigenvalue (w.r.t. $\omega_{\rm FS}$) bounded below by $1$.
We observe that this single positive Levi direction
suffices to make the weighted metric uniformly
Nakano positive (see Lemma \ref{lem-projections}
with $p=r=1$ and $a=1$).
Demailly's theorem, together with the boundedness of
the weight, then yields the desired vanishing and
gives an affirmative answer to
Question \ref{ques-shaw-intro}.

More generally, the same argument applies to domains in
$\mathbb P^n$ admitting a complete K\"ahler metric and
a bounded smooth plurisubharmonic weight with at least
$r$ uniformly positive Levi eigenvalues (w.r.t. $\omega_{\rm FS}$).

\begin{theorem}\label{thm-weight}
Let $\Omega\subset\mathbb P^n$ be a domain, where $n\geq2$.
Fix $1\leq r\leq n$ and $a>0$. Suppose that
\begin{enumerate}[\rm{(}1\rm{)}]
\item $\Omega$ admits a complete K\"ahler metric;
\item there is a bounded smooth plurisubharmonic function
$\varphi$ on $\Omega$ with the following property:
for each $x\in\Omega$, let
$\lambda_1(x)\geq\cdots\geq\lambda_n(x)\geq0$
be the eigenvalues   of the
Levi form $\sqrt{-1}\partial\bar\partial\varphi$
with respect to $\omega_{\rm FS}$ at $x$.
Then
$\lambda_j(x)\geq a$ for $1\leq j\leq r$
and each $x\in\Omega$.
\end{enumerate}
Then, for $1\leq p\leq r$, $1\leq q\leq n$ and each
$f\in\ker(\bar\partial^{p,q})$, there exists
$u\in\operatorname{Dom}(\bar\partial^{p,q-1})$
such that $\bar\partial u=f$ and
\begin{equation}\label{eq-main-estimate}
 \|u\|^2
 \leq
 \frac{n+1-p+a}{q\,a(r-p+1)}
 \exp\bigl(\sup_{\Omega}\varphi-\inf_{\Omega}\varphi)\,
 \|f\|^2,
\end{equation}
 where the $L^2$-norm $||\bullet||$ is with respect to the restriction of the Fubini--Study metric to $\Omega$.
That is to say, 
$H^{p,q}_{(2)}(\Omega)=0$ for $1\leq p\leq r$
and $1\leq q\leq n$.  
\end{theorem}

Note that a pseudoconvex  domain in $\mathbb C^n$ is Stein and therefore
admits a complete K\"ahler metric  (\cite[Theorem 1.3]{Dem82}).  We then obtain the following corollary of Theorem \ref{thm-weight}.

\begin{corollary}\label{cor-coordinates}
Let $\Omega$ be a pseudoconvex domain in a standard affine chart
$\mathbb C^{n-r}_z\times\mathbb C^r_w\subset\mathbb P^n$,
where $1\leq r\leq n$ and
\eqref{eq-fubini-study} holds with $\zeta=(z,w)$.
Suppose that $|w|<R$ on $\Omega$.
Then, for $1\leq p\leq r$, $1\leq q\leq n$ and each
$f\in\ker(\bar\partial^{p,q})$, there exists
$u\in\operatorname{Dom}(\bar\partial^{p,q-1})$
such that $\bar\partial u=f$ and
\[
 \|u\|^2
 \leq
 \frac{n+2-p}{q(r-p+1)}
 e^{R^2}\|f\|^2.
\]
Consequently,
$H^{p,q}_{(2)}(\Omega)=0$ for $1\leq p\leq r$
and $1\leq q\leq n$.
\end{corollary}

We can  apply Corollary \ref{cor-coordinates} to obtain an
affirmative answer to Question  \ref{ques-shaw-intro}.

\begin{corollary}[{Answer to Question \ref{ques-shaw-intro}}]\label{cor-hartogs}
Let $\mathbb H_n^+$ be the  projective Hartogs triangles in
$\mathbb P^n$ with
 $n\geq2$.
Then, for $1\leq q\leq n$ and each
$f\in\ker(\bar\partial^{1,q})$, there exists
$u\in\operatorname{Dom}(\bar\partial^{1,q-1})$
such that $\bar\partial u=f$ and
\[
 \|u\|^2
 \leq
 \frac{(n+1)e}{q}\|f\|^2.
\]
Consequently,
\[
 H^{1,q}_{(2)}(\mathbb H_n^+)=0
 \qquad\text{for }1\leq q\leq n.
\]
In particular, Question \ref{ques-shaw-intro} has an affirmative answer.
\end{corollary}

Note that Shaw proved that
$H^{2,1}_{(2)}(\mathbb H_2^\pm)$ is infinite dimensional
(\cite[Lemma 4.3(3)]{Sh25}).  This suggests a more general
nonvanishing phenomenon and, in particular, yields the question
whether the range of holomorphic degrees $p\leq r$ in
Theorem \ref{thm-weight} and Corollary \ref{cor-coordinates}
is optimal.  Motivated by the phenomenon that Bergman spaces
of holomorphic forms can describe  obstructions to
solving $\bar\partial$-equations in $L^2$ and Sobolev spaces
(e.g., \cite[Theorem 2.1]{Sh11},
\cite[Theorem 1.1]{Ch24}), we prove  that
$H^{r+1,1}_{(2)}
\bigl(\mathbb C^{n-r}\times\mathbb D^r\bigr)$
is infinite dimensional (w.r.t. the restricted
Fubini--Study metric), via constructing an 
injective linear map from the Bergman space
$A^2(\mathbb D^r)$ into this $L^2$ Dolbeault cohomology group.
This shows that the range $p\leq r$ in
Theorem \ref{thm-weight} and Corollary \ref{cor-coordinates}
is optimal.

\begin{theorem}\label{thm-sharp}
Let $1\leq r<n$, and view
$\mathbb C^{n-r}_z\times\mathbb D^r_w$
as a domain in the standard affine chart of $\mathbb P^n$, where $\mathbb D^r_w$ denotes the  unit polydisc in $\mathbb C^{r}$.
Set
\[
 \omega_{n-r}
 =
 \sqrt{-1}\partial\bar\partial
 \log(1+|z|^2),
\]
where  $A^2(\mathbb D^r)$ denote
the Bergman space  (w.r.t. the Euclidean metric).
Then there is an injective linear map
\begin{equation}\label{eq-bergman-map}
 \begin{aligned}
 A^2(\mathbb D^r)&\longrightarrow
 H^{r+1,1}_{(2)}
 \bigl(\mathbb C^{n-r}_z\times\mathbb D^r_w\bigr),\\
 h&\longmapsto
 [\,h(w)\omega_{n-r}(z)\wedge
 dw_1\wedge\cdots\wedge dw_r\,].
 \end{aligned}
\end{equation}
In particular,
$H^{r+1,1}_{(2)}
 \bigl(\mathbb C^{n-r}\times\mathbb D^r\bigr)$ (w.r.t. the restricted Fubini–Study metric)
is infinite dimensional.
\end{theorem}

\begin{remark}\label{rem-product-metric}
Note that the $L^2$ K\"unneth formula of
Chakrabarti--Shaw
\cite[Theorem 1.1 and Corollary 6.5]{CS11} can imply that $H^{r+1,1}_{(2)}
\bigl(
\mathbb C^{n-r}\times\mathbb D^r;
g_0
\bigr)$
is infinite dimensional, where $g_0=
g_{{\rm FS},n-r}|_{\mathbb C^{n-r}}
\oplus g_{{\rm E},r}$
is the product of the restriction of the Fubini--Study metric
from $\mathbb P^{n-r}$ to $\mathbb C^{n-r}$ and the Euclidean
metric on $\mathbb D^r$. However
in Theorem \ref{thm-sharp},
$\mathbb C^{n-r}\times\mathbb D^r$ is viewed as a domain in
$\mathbb P^n$ and is equipped with the Fubini--Study metric
induced from the ambient projective space.
Thus the above K\"unneth formula does not directly imply
Theorem \ref{thm-sharp}.
\end{remark}

The rest of the paper is organized as follows.
In Section \ref{sec-curvature}, we recall Demailly's
$L^2$ existence theorem, establish the required curvature
estimates, and prove Theorem \ref{thm-weight} and
Corollaries \ref{cor-coordinates} and \ref{cor-hartogs}.
In Section \ref{sec-optimality}, we prove
Theorem \ref{thm-sharp} by constructing an explicit
injection from a Bergman space into an $L^2$ Dolbeault
cohomology group, thus showing that the condition $p\leq r$   in Theorem \ref{thm-weight} and
Corollary \ref{cor-coordinates} is indispensable.

\section{Curvature  estimate and proofs of main results}\label{sec-curvature}

\subsection{Curvature operators and Demailly's $L^2$ existence theorem}
\label{subsec-demailly-estimate}
We refer the reader to \cite[Chapters VII and VIII]{Dem12} for a systematic treatment of the material in this subsection.

Let $(X,\omega)$ be a Hermitian manifold of dimension $n$,
and let $(E,h)$ be a Hermitian holomorphic vector bundle
of rank $\rho$.
Write $\Theta(E,h)$ for the
$\operatorname{End}(E)$-valued Chern curvature
$(1,1)$-form of $(E,h)$.
All norms, inner products and $L^2$ spaces are computed
with respect to $\omega$ and $h$, and
$dV_\omega=\omega^n/n!$.

For an integer $m\geq1$, the bundle $(E,h)$ is called
\emph{$m$-semipositive} (e.g., 
\cite[Chapter VII, Definition 6.5]{Dem12}) if, at each $x\in X$,
\[
 \sum_{j,k=1}^n
 \left\langle
 \Theta_{j\bar k}(E,h)v_j,v_k
 \right\rangle_h
 \geq0
\]
for each tuple $v_1,\ldots,v_n\in E_x$ satisfying
$\dim_{\mathbb C}\operatorname{span}
\{v_1,\ldots,v_n\}\leq m$,
where
$$\Theta_{j\bar k}(E,h)
=\Theta(E,h)_x(e_j,\overline{e_k})$$
in an $\omega$-unitary tangent frame
$e_1,\ldots,e_n$.   When $m=\min\{n,\rho\}$, this condition is precisely
Nakano semipositivity. 

Let $L_\omega$ denote the \emph{Lefschetz operator},
defined by $L_\omega u=\omega\wedge u$,
and let $\Lambda_\omega$ denote its pointwise
 adjoint, called the \emph{dual Lefschetz operator}.
On the space of  $E$-valued $(n,q)$-forms with $q\geq 1$, we define the pointwise Hermitian  operator $A_q=A_q(E,h,\omega)$ by
\begin{equation}\label{eq-curvature-operator}
 A_q u
 =
 \sqrt{-1}\Theta(E,h)\wedge(\Lambda_\omega u).
\end{equation}
For $E$-valued $(n,1)$-forms $ u=\sum_{i=1}^n
 \sigma\wedge\overline{\theta^i}\otimes u_i,$
we have
\begin{equation}\label{A1u-u-Theta}
    \langle A_1u,u\rangle_{\omega,h}
 =
 \sum_{i,j=1}^n
 \left\langle
 \Theta_{i\bar j}(E,h)u_i,u_j
 \right\rangle_h.
\end{equation}
When acting on $E$-valued $(n,q)$-forms, we have
\begin{equation}\label{Lomega-Lamb-qI}
     L_\omega\Lambda_\omega=qI.
\end{equation}
Moreover, If $(E,h)$ is $m$-semipositive and
$m\geq\min\{n-q+1,\rho\}$, then $A_q\geq0$ (\cite[Chapter VII, (7.1)]{Dem12}).

For $A_q\geq0$, the expression
$\langle A_q^{-1}g,g\rangle_{\omega,h}$ at $x$
denotes the smallest constant $M\geq0$ such that
\[
 |\langle g,v\rangle_{\omega,h}|^2
 \leq
 M\langle A_qv,v\rangle_{\omega,h}
\]
for each $v\in\Lambda^{n,q}T_x^*X\otimes E_x$
(\cite[Chapter VIII, Section 4, before (4.4)]{Dem12}).
If no finite constant has this property, the expression
is defined to be $+\infty$
(\cite[pp. 468--469, (3.1)]{Dem82}).
When $A_q$ is positive definite, this agrees with
the quadratic form defined by the usual inverse.

The following $L^2$ existence  theorem of Demailly  plays an important role in proving Theorem  \ref{thm-weight}.

\begin{lemma}[{\cite[Theorem 4.1 and Remark 4.2]{Dem82} or \cite[Chapter VIII, Theorem 6.1]{Dem12}}]\label{lem-demailly-existence}
Let $(X,\widehat\omega)$ be a complete K\"ahler
manifold of dimension $n$, and let $\omega$ be
another K\"ahler metric on $X$, which need not
be complete.
Let $(E,h)$ be a Hermitian holomorphic vector bundle
of rank $\rho$ that is $m$-semipositive.
Fix $1\leq q\leq n$ such that
$m\geq\min\{n-q+1,\rho\}$.  Then for each $g\in L^2_{n,q}(X,E;\omega,h)$ satisfying
$\bar\partial g=0$ and
\[
 \int_X
 \langle A_q^{-1}g,g\rangle_{\omega,h}
 \,dV_\omega<+\infty,
\]
there exists $f\in L^2_{n,q-1}(X,E;\omega,h)$
such that $\bar\partial f=g$ and
\begin{equation}\label{eq-demailly-estimate}
 \|f\|_{\omega,h}^2
 \leq
 \int_X
 \langle A_q^{-1}g,g\rangle_{\omega,h}
 \,dV_\omega.
\end{equation}
\end{lemma}

\subsection{Curvature estimate}
\label{subsec-induced-curvature}

We now work on $\mathbb P^n$ with the
Fubini--Study metric $\omega_{\rm FS}$.
For $1\leq p\leq n$, let
$$G_p=K_{\mathbb P^n}^{-1}\otimes\Omega_{\mathbb P^n}^p$$
and denote its induced Hermitian metric by $h_p$.
For bundles with their induced metrics, we omit
the metric from the curvature notation and write
$\Theta(E)$ and $\Theta_{i\bar j}(E)$.
We also write
$A_1(G_p)=A_1(G_p,h_p,\omega_{\rm FS})$. We make the convention that 
pointwise norms $|\bullet|$ without subscripts in this subsection
are computed using the induced metrics.

Via the canonical isomorphism (isometry)
\begin{equation}\label{eq-bundle-identification}
 \Omega_{\mathbb P^n}^p
 \cong K_{\mathbb P^n}\otimes G_p,
\end{equation}
we can consider $(p,q)$-forms as $G_p$-valued
$(n,q)$-forms.

Fix $x\in\mathbb P^n$ and an $\omega_{\rm FS}$-unitary
tangent frame $e_1,\ldots,e_n$, with the dual coframe
$\theta^1,\ldots,\theta^n$.
Write $e_j\lrcorner\alpha$ for contraction of a form
$\alpha$ with $e_j$.
For a tuple $U=(U_1,\ldots,U_n)$ with
$U_j\in\Lambda^p((T_x^{1,0}\mathbb P^n)^*)$,
define
$CU=\sum_{j=1}^n(e_j\lrcorner U_j)$,
and let $C^*$ denote its Hermitian adjoint.
Since the frame is unitary,  we obtain
\begin{equation}\label{Cstar-alpha}
    C^*\alpha
 =
(\theta^1\wedge\alpha,\ldots,\theta^n\wedge\alpha)
\end{equation}
for each
$\alpha\in\Lambda^{p-1}((T_x^{1,0}\mathbb P^n)^*)$.
The exterior algebra identity
$\sum_{j=1}^n\theta^j\wedge(e_j\lrcorner\alpha)
=(p-1)\alpha$ then gives that
\[
 CC^*\alpha
 =
 \sum_{j=1}^n e_j\lrcorner(\theta^j\wedge\alpha)
 =
 (n+1-p)\alpha.
\]
Thus $CC^*=(n+1-p)I$.


Put $\sigma=\theta^1\wedge\cdots\wedge\theta^n$,
and let $\sigma^{-1}\in(K_{\mathbb P^n}^{-1})_x$
be its dual vector such that
$\sigma^{-1}(\sigma)=1$.
Then any $G_p$-valued $(n,1)$-form at $x$ can be
written uniquely as
\[
 \sum_{j=1}^n
 \sigma\wedge\overline{\theta^j}
 \otimes(\sigma^{-1}\otimes U_j),
\]
where
$U=(U_1,\ldots,U_n)$ and
$U_j\in\Lambda^p((T_x^{1,0}\mathbb P^n)^*)$.
Clearly, the above correspondence gives an isometric identification
between the space of such tuples $U$ and
$\Lambda^{n,1}T_x^*\mathbb P^n\otimes(G_p)_x$.
We use the same symbol $U$ for a tuple and its corresponding
$G_p$-valued $(n,1)$-form. Then   we have
\[
 |U|^2=\sum_{j=1}^n|U_j|^2.
\]

We first express the curvature operator $A_1(G_p)$ in terms of the contraction map $C$ and establish the Nakano semipositivity of $G_p$.

\begin{lemma}\label{lem-curvature}
Under the above isometric identification, the curvature operator satisfies
\begin{equation}\label{A1GpU-Cstar}
 A_1(G_p)U=(n+1-p)U-C^*CU
\end{equation}
for any $U=(U_1,\ldots,U_n)$ with
$U_j\in\Lambda^p((T_x^{1,0}\mathbb P^n)^*)$.
Consequently,
\begin{equation}\label{eq-curvature-form}
 \langle A_1(G_p)U,U\rangle
 =
 (n+1-p)|U|^2-|CU|^2 \geq 0.
\end{equation}
\end{lemma}

\begin{proof}
Write $T=T^{1,0}\mathbb P^n$.
With the normalization \eqref{eq-fubini-study},
the Fubini--Study curvature formula
(e.g., \cite[(2.1)--(2.2) and (7.5)]{LSY13}) gives
\[
 \left\langle
 \Theta_{i\bar j}(T)e_k,e_\ell
 \right\rangle
 =
 R_{i\bar j k\bar\ell}
 =
 \delta_{ij}\delta_{k\ell}
 +
 \delta_{i\ell}\delta_{kj}.
\]
Since this frame is orthonormal,
\[
\begin{aligned}
 \Theta_{i\bar j}(T)e_k
 &=
 \sum_{\ell=1}^n
 \left\langle
 \Theta_{i\bar j}(T)e_k,e_\ell
 \right\rangle e_\ell\\
 &=
 \sum_{\ell=1}^n
 \bigl(
 \delta_{ij}\delta_{k\ell}
 +
 \delta_{i\ell}\delta_{kj}
 \bigr)e_\ell\\
 &=
 \delta_{ij}e_k+\delta_{kj}e_i.
\end{aligned}
\]

Since the coframe is unitary,
$|\sigma|=|\sigma^{-1}|=1$.
Thus the induced tensor product metrics satisfy
\[
 \left\langle
 \sigma^{-1}\otimes\alpha,
 \sigma^{-1}\otimes\beta
 \right\rangle
 =
 \langle\alpha,\beta\rangle
\]
for any $(p,0)$-forms $\alpha$ and $\beta$ at $x$.
The curvature formulas for duals, exterior powers
and determinants
(e.g., \cite[Chapter V, $(4.3')$, $(4.5')$ and (4.6)]{Dem12})
give
\[
\begin{aligned}
 \Theta_{i\bar j}(\Omega_{\mathbb P^n}^p)\alpha
 &=
 -p\delta_{ij}\alpha
 -\theta^j\wedge(e_i\lrcorner\alpha),\\
 \Theta_{i\bar j}(K_{\mathbb P^n}^{-1})\sigma^{-1}
 &=
 (n+1)\delta_{ij}\sigma^{-1},
\end{aligned}
\]
where the second equality follows from the facts that
$K_{\mathbb P^n}^{-1}=\det T$ and
$\operatorname{tr}(\Theta_{i\bar j}(T))
=(n+1)\delta_{ij}$.
Applying the tensor product curvature formula
(e.g., \cite[Chapter V, $(4.2')$]{Dem12}), we obtain
\[
\begin{aligned}
 &\Theta_{i\bar j}(G_p)
   (\sigma^{-1}\otimes\alpha)\\
 &\qquad=
 \bigl(
 \Theta_{i\bar j}(K_{\mathbb P^n}^{-1})\sigma^{-1}
 \bigr)\otimes\alpha
 +
 \sigma^{-1}\otimes
 \bigl(
 \Theta_{i\bar j}(\Omega_{\mathbb P^n}^p)\alpha
 \bigr)\\
 &\qquad=
 (n+1-p)\delta_{ij}
   (\sigma^{-1}\otimes\alpha)
 -
 \sigma^{-1}\otimes
 \bigl(\theta^j\wedge(e_i\lrcorner\alpha)\bigr).
\end{aligned}
\]
Taking $\alpha=U_i$, we obtain
\begin{equation}\label{Thetaij-sigmainver-Ui}
 \Theta_{i\bar j}(G_p)
   (\sigma^{-1}\otimes U_i)
 =
 (n+1-p)\delta_{ij}
   (\sigma^{-1}\otimes U_i)
 -
 \sigma^{-1}\otimes
 \bigl(\theta^j\wedge(e_i\lrcorner U_i)\bigr).
\end{equation}

For any $G_p$-valued $(n,1)$-form
\[
 U=
 \sum_{i=1}^n
 \sigma\wedge\overline{\theta^i}
 \otimes(\sigma^{-1}\otimes U_i),
\]
the definition \eqref{eq-curvature-operator} shows that
the coefficient of
$\sigma\wedge\overline{\theta^j}$ in $A_1(G_p)U$
is
$\sum_{i=1}^n
 \Theta_{i\bar j}(G_p)
 (\sigma^{-1}\otimes U_i)$.
Using the curvature formula
\eqref{Thetaij-sigmainver-Ui} and the definition
$CU=\sum_{i=1}^n e_i\lrcorner U_i$, we obtain
\[
\begin{aligned}
 &\sum_{i=1}^n
 \Theta_{i\bar j}(G_p)
   (\sigma^{-1}\otimes U_i)\\
 &\qquad=
 \sigma^{-1}\otimes
 \bigl((n+1-p)U_j-\theta^j\wedge CU\bigr)\\
 &\qquad=
 \sigma^{-1}\otimes
 \left[
 \bigl((n+1-p)I-C^*C\bigr)U
 \right]_j,
\end{aligned}
\]
where the last equality follows from
$[C^*CU]_j=\theta^j\wedge CU$.
Thus, under the above isometric identification,
$A_1(G_p)U=(n+1-p)U-C^*CU$,
which proves \eqref{A1GpU-Cstar}.

Taking the inner product of \eqref{A1GpU-Cstar}
with $U$ and using the definition of the Hermitian
adjoint $C^*$, we obtain
\[
\begin{aligned}
 \langle A_1(G_p)U,U\rangle
 &=
 (n+1-p)|U|^2-\langle C^*CU,U\rangle\\
 &=
 (n+1-p)|U|^2-|CU|^2.
\end{aligned}
\]
This proves the first part of \eqref{eq-curvature-form}.

For each
$\gamma\in\Lambda^{p-1}((T_x^{1,0}\mathbb P^n)^*)$,
the identity $CC^*=(n+1-p)I$ gives
$|C^*\gamma|^2=(n+1-p)|\gamma|^2$.
By duality and the Cauchy--Schwarz inequality,
\[
 |CU|^2
 =
 \sup_{|\gamma|=1}
 |\langle U,C^*\gamma\rangle|^2
 \leq(n+1-p)|U|^2.
\]
Thus $ \langle A_1(G_p)U,U\rangle\geq 0$ for  any $U=(U_1,\ldots,U_n)$ with
$U_j\in\Lambda^p((T_x^{1,0}\mathbb P^n)^*)$. This proves the second part of  \eqref{eq-curvature-form}.
\end{proof}

We next  show that for $p\leq r$, adding
uniform positivity in just $r$ directions
suffices to make the nonnegative quadratic form in
\eqref{eq-curvature-form} uniformly positive definite.  In particular, a plurisubharmonic weight with $r$ uniformly positive Levi directions will yield a uniform positive lower bound for the weighted  curvature.

\begin{lemma}\label{lem-projections}
With the notation introduced before Lemma \ref{lem-curvature},
let $1\leq p\leq r\leq n$ be integers and let $a>0$.
Then for each tuple $U=(U_1,\ldots,U_n)$ with
$U_j\in\Lambda^p((T_x^{1,0}\mathbb P^n)^*)$, we have
\begin{equation}\label{eq-weighted-curvature}
 (n+1-p)|U|^2-|CU|^2
 +a\sum_{j=1}^r|U_j|^2
 \geq
 \frac{a(r-p+1)}{n+1-p+a}|U|^2.
\end{equation}
\end{lemma}

\begin{proof}
For each
$\alpha\in\Lambda^{p-1}((T_x^{1,0}\mathbb P^n)^*)$,
we have
\[
 \sum_{j=1}^n|\theta^j\wedge\alpha|^2
 =(n+1-p)|\alpha|^2
\]
and
\begin{equation}\label{eq-kernel-control}
 \begin{aligned}
 \sum_{j=1}^r|\theta^j\wedge\alpha|^2
 &=r|\alpha|^2-
   \sum_{j=1}^r|e_j\lrcorner\alpha|^2\\
 &\geq(r-p+1)|\alpha|^2,
 \end{aligned}
\end{equation}
where the last inequality follows from
$\sum_{j=1}^n|e_j\lrcorner\alpha|^2
=(p-1)|\alpha|^2$.

Note that 
$$|\langle CU,\alpha\rangle|^2
 =
 \left|
 \sum_{j=1}^n
 \langle U_j,\theta^j\wedge\alpha\rangle
 \right|^2.$$
We then apply the Cauchy--Schwarz inequality in the form
\[
 \left|
 \sum_{j=1}^n
 \left\langle
 \sqrt{b_j}\,x_j,
 \frac{1}{\sqrt{b_j}}\,y_j
 \right\rangle
 \right|^2
 \leq
 \left(
 \sum_{j=1}^n b_j|x_j|^2
 \right)
 \left(
 \sum_{j=1}^n\frac{|y_j|^2}{b_j}
 \right),
\]
to $ \left|
 \sum_{j=1}^n
 \langle U_j,\theta^j\wedge\alpha\rangle
 \right|^2$ with
\[
 b_j=
 \begin{cases}
 n+1-p+a,&1\leq j\leq r,\\
 n+1-p,&r<j\leq n,
 \end{cases}
\]
and  obtain

\[
\begin{aligned}
 |\langle CU,\alpha\rangle|^2
 &\leq
 \left(
 (n+1-p)|U|^2
 +a\sum_{j=1}^r|U_j|^2
 \right)\\
 &\qquad\times
 \left(
 1-
 \frac{a(r-p+1)}
 {(n+1-p)(n+1-p+a)}
 \right)|\alpha|^2.
\end{aligned}
\]
Taking the supremum over $|\alpha|=1$, we obtain
\[
 |CU|^2
 \leq
 \left(
 1-
 \frac{a(r-p+1)}
 {(n+1-p)(n+1-p+a)}
 \right)
 \left(
 (n+1-p)|U|^2
 +a\sum_{j=1}^r|U_j|^2
 \right),
\]
which gives \eqref{eq-weighted-curvature}.
\end{proof}

\subsection{Proofs of main results}\label{sec-vanishing}

Now we are ready to prove the main results.

\begin{proof}[Proof of Theorem \ref{thm-weight}]
Fix $1\leq p\leq r$, $1\leq q\leq n$ and
$f\in\ker(\bar\partial^{p,q})$.
Restrict the bundle $G_p$ to $\Omega$ and equip it with the weighted metric
$$h_{p,\varphi}=e^{-\varphi}h_p,$$ and write
$$A_{q,\varphi}
=A_q(G_p,h_{p,\varphi},\omega_{\rm FS}).$$
Then the curvature of this metric is
$$\Theta(G_p,h_{p,\varphi})
=\Theta(G_p,h_p)+\partial\bar\partial\varphi\otimes I.$$

At an arbitrary point $x\in\Omega$, choose an
$\omega_{\rm FS}$-unitary tangent frame diagonalizing
the Levi form of $\varphi$.
Let
$\lambda_1\geq\cdots\geq\lambda_n\geq0$
be the corresponding Levi eigenvalues.
By assumption,
$\lambda_1\geq\cdots\geq\lambda_r\geq a$.
In this frame,
$\Theta_{i\bar j}(G_p,h_{p,\varphi})
=\Theta_{i\bar j}(G_p,h_p)+\lambda_i\delta_{ij}I$.
We consider any $G_p$-valued
$(n,1)$-form at $x$ as a tuple $U=(U_1,\ldots,U_n)$. 
In the following pointwise calculation, unsubscripted
norms are computed using the original induced metrics.

It follows from 
\eqref{eq-curvature-operator},
\eqref{eq-curvature-form} and
\eqref{eq-weighted-curvature} that
\[
\begin{aligned}
 &\langle A_{1,\varphi}U,U
   \rangle_{\omega_{\rm FS},h_p}\\
 &\qquad=
 (n+1-p)|U|^2-|CU|^2
 +\sum_{j=1}^n\lambda_j|U_j|^2\\
 &\qquad\geq
 (n+1-p)|U|^2-|CU|^2
 +a\sum_{j=1}^r|U_j|^2\\
 &\qquad\geq
 c|U|^2,
\end{aligned}
\]
where
\begin{equation}\label{eq-curvature-constant}
 c=\frac{a(r-p+1)}{n+1-p+a}>0.
\end{equation}
Consequently, we obtain
\[
\begin{aligned}
 &\langle A_{1,\varphi}U,U
   \rangle_{\omega_{\rm FS},h_{p,\varphi}}\\
 &\qquad=
 e^{-\varphi(x)}
 \langle A_{1,\varphi}U,U
   \rangle_{\omega_{\rm FS},h_p}\\
 &\qquad\geq
 c\,e^{-\varphi(x)}|U|^2
 =
 c\,|U|_{\omega_{\rm FS},h_{p,\varphi}}^2.
\end{aligned}
\]
Since $x$ and $U$ are arbitrary, we obtain
\begin{equation}\label{eq-weighted-operator}
 A_{1,\varphi}\geq cI,
 \qquad
 A_{1,\varphi}^{-1}\leq c^{-1}I.
\end{equation}
By \eqref{A1u-u-Theta}, this  gives 
uniform Nakano positivity of $(G_p,h_{p,\varphi})$. In particular, $(G_p,h_{p,\varphi})$ is  $m$-semipositive for
$m=\min\{n,\operatorname{rank}G_p\}$.

The pointwise Hermitian form
$\sqrt{-1}\Theta(G_p,h_{p,\varphi})
-c\omega_{\rm FS}\otimes I$
is Nakano semipositive by \eqref{eq-weighted-operator}.
By the pointwise  curvature formula  on $G_p$-valued
$(n,q)$-forms, its induced operator
$A_{q,\varphi}-cL_{\omega_{\rm FS}}\Lambda_{\omega_{\rm FS}}$
is therefore semipositive.  By \eqref{Lomega-Lamb-qI},  $L_{\omega_{\rm FS}}\Lambda_{\omega_{\rm FS}}=qI$ on these forms.  Then we obtain
\begin{equation}\label{eq-higher-curvature}
 A_{q,\varphi}\geq qcI,
 \qquad
 A_{q,\varphi}^{-1}\leq\frac{1}{qc}I.
\end{equation}

Since $\rho=\operatorname{rank}G_p=\binom{n}{p}$
and the weighted bundle is Nakano positive,
it is $m$-semipositive for
$m=\min\{n,\rho\}$, which satisfies
$m\geq\min\{n-q+1,\rho\}$. Then we can apply
Lemma \ref{lem-demailly-existence}  to the bundle $E=G_p|_\Omega$ with the weighted metric $h_{p,\varphi}$. 
Again via the isometry \eqref{eq-bundle-identification},
we regard $f$ as a $G_p$-valued $(n,q)$-form.
The boundedness of $\varphi$ gives the integrability $\int_\Omega |f|^2e^{-\varphi}\,dV_{\rm FS}
 <+\infty,$
and thus we have
\[
\begin{aligned}
 &\int_\Omega
 \left\langle
 A_{q,\varphi}^{-1}f,f
 \right\rangle_{\omega_{\rm FS},h_{p,\varphi}}
 \,dV_{\rm FS}\\
 &\qquad \leq
 \frac{1}{qc}\int_\Omega |f|^2e^{-\varphi}\,dV_{\rm FS}
 <+\infty.
\end{aligned}
\]
Then 
Lemma \ref{lem-demailly-existence} therefore gives
$u\in L^2_{n,q-1}(\Omega,G_p;\omega_{\rm FS},h_{p,\varphi})$
such that $\bar\partial u=f$ and
\[
 \|u\|_{\omega_{\rm FS},h_{p,\varphi}}^2
 \leq
 \frac{1}{qc}\int_\Omega |f|^2e^{-\varphi}\,dV_{\rm FS}.
\]

Since $\varphi$ is bounded, we have
\[
 e^{-\sup_\Omega\varphi}
 \leq e^{-\varphi}
 \leq e^{-\inf_\Omega\varphi}
\]
on $\Omega$.
Therefore,
\[
\begin{aligned}
 \|u\|^2
 &=
 \int_\Omega |u|^2\,dV_{\rm FS}\\
 &\leq
 e^{\sup_\Omega\varphi}
 \int_\Omega |u|^2e^{-\varphi}\,dV_{\rm FS}\\
 &\leq
 \frac{1}{qc}
 e^{\sup_\Omega\varphi-\inf_\Omega\varphi}
 \|f\|^2.
\end{aligned}
\]
Using \eqref{eq-curvature-constant}, this proves
\eqref{eq-main-estimate}.

As $f\in\ker(\bar\partial^{p,q})$ was arbitrary, we obtain
$\operatorname{im}(\bar\partial^{p,q-1})
=\ker(\bar\partial^{p,q})$ and thus
$H^{p,q}_{(2)}(\Omega)=0$.  This completes the proof of Theorem \ref{thm-weight}. 
\end{proof}

\begin{proof}[Proof of Corollary \ref{cor-coordinates}]
Set $\varphi(z,w)=|w|^2$ and put
$\rho=1+|\zeta|^2=1+|z|^2+|w|^2$.
With respect to the Fubini--Study metric, the nonzero
Levi eigenvalues of $\varphi$ are
\[
 \rho(1+|w|^2),\qquad
 \underbrace{\rho,\ldots,\rho}_{r-1}.
\]
Thus, 
$\lambda_1\geq\cdots\geq\lambda_r\geq 1$
are the first $r$ Levi eigenvalues of $\varphi$ with respect to
$\omega_{\rm FS}$.
Moreover,
$\sup_\Omega\varphi-\inf_\Omega\varphi \leq R^2$.
Note that the pseudoconvex affine domain $\Omega$ admits a complete K\"ahler metric (e.g., \cite[Chapter VIII, Theorem 5.2]{Dem12}).  Then 
Theorem \ref{thm-weight}  gives
the assertion.
\end{proof}

\begin{proof}[Proof of Corollary \ref{cor-hartogs}]
The defining inequality \eqref{def-intro-tri} for $\mathbb H_n^+$ implies
$Z_n\neq0$. In the standard affine chart
$\{Z_n\neq0\}$, put
\[
 z=
 \left(
 \frac{Z_0}{Z_n},\ldots,\frac{Z_{n-2}}{Z_n}
 \right)
 \in\mathbb C^{n-1},
 \qquad
 w=\frac{Z_{n-1}}{Z_n}\in\mathbb C.
\]
Then
$\mathbb H_n^+
 \cong
 \mathbb C^{n-1}_z\times\mathbb D_w$ is pseudoconvex, and
Corollary \ref{cor-coordinates} applies with
$r=1$, $p=1$ and $R=1$. Thus, for
$1\leq q\leq n$ and each
$f\in\ker(\bar\partial^{1,q})$, there exists
$u\in\operatorname{Dom}(\bar\partial^{1,q-1})$
such that $\bar\partial u=f$ and
\[
 \|u\|^2
 \leq
 \frac{(n+1)e}{q}\|f\|^2.
\]
Consequently,
\[
 H^{1,q}_{(2)}(\mathbb H_n^+)=0
 \qquad\text{for }1\leq q\leq n.
\]
\end{proof}

\section{Injectivity from Bergman space to  $L^2$ Dolbeault cohomology groups }\label{sec-optimality}
In this section, we construct in Theorem \ref{thm-sharp} an explicit injective linear map
from a Bergman space into an $L^2$ Dolbeault cohomology group.
This yields the infinite dimensionality of certain $L^2$
Dolbeault cohomology groups and shows that the condition $p\leq r$   in Theorem \ref{thm-weight} and
Corollary \ref{cor-coordinates} is indispensable. 

The proof of Theorem \ref{thm-sharp} is based on the following two observations.

\begin{observation}\label{lem-affine-obstruction}
For $n-r\geq1$, there is no
$\beta\in L^2_{1,0}(\mathbb C^{n-r},\omega_{n-r})$
satisfying $\bar\partial\beta=\omega_{n-r}$
in the distributional sense.
\end{observation}

\begin{proof}
Assume by contradiction that such a $\beta$ exists.
Let
$H_\infty=\mathbb P^{n-r}\setminus\mathbb C^{n-r}$
be the hyperplane at infinity.
The form $\omega_{n-r}$ is the restriction of the
Fubini--Study form on $\mathbb P^{n-r}$.
By the $L^2$ removability theorem
(\cite[Lemma 6.9]{Dem82}), the form $\beta$ can extend
across $H_\infty$ as an $L^2$ form satisfying
$\bar\partial\beta=\omega_{n-r}$ on $\mathbb P^{n-r}$
in the distributional sense.
Since $\bar\partial\omega_{n-r}=0$, we obtain
$\bar\partial
 \bigl(\beta\wedge\omega_{n-r}^{n-r-1}\bigr)
 =
 \omega_{n-r}^{n-r}$
on $\mathbb P^{n-r}$, contradicting the positivity of the
class $[\omega_{n-r}^{n-r}]$.
\end{proof}

\begin{observation}\label{lem-metric-comparison}
The restriction of the Fubini--Study metric
$g_{{\rm FS},n}$ on $\mathbb P^n$ to
$\mathbb C^{n-r}_z\times\mathbb D^r_w
\subseteq \mathbb C^n\subseteq \mathbb P^n$ is uniformly
equivalent to the diagonal type  Hermitian metric
\begin{equation}\label{eq-product-metric}
 \widetilde g
 =
 g_{{\rm FS},n-r}(z)\oplus
 (1+|z|^2)^{-1}g_{{\rm E},r}(w),
\end{equation}
where $g_{{\rm E},r}$ denotes the standard Euclidean
 metric on $\mathbb C^r$.
\end{observation}

\begin{proof}
For a cotangent vector
\[
 \xi
 =
 \sum_{j=1}^{n-r} b_j\,dz_j
 +
 \sum_{\alpha=1}^{r} d_\alpha\,dw_\alpha,
\]
with
$b=(b_1,\ldots,b_{n-r})$ and
$d=(d_1,\ldots,d_r)$, the squared dual norms are
\[
\begin{aligned}
 |(b,d)|_{\rm FS}^2
 &=
 (1+|z|^2+|w|^2)
 \bigl(
 |b|^2+|d|^2+|z\cdot b+w\cdot d|^2
 \bigr),\\
 |(b,d)|_{\widetilde g}^2
 &=
 (1+|z|^2)
 \bigl(
 |b|^2+|d|^2+|z\cdot b|^2
 \bigr).
\end{aligned}
\]
Since $|w|^2<r$, we have
\[
 1
 \leq
 \frac{1+|z|^2+|w|^2}{1+|z|^2}
 \leq
 1+r.
\]
Using
$|s+t|^2\leq2|s|^2+2|t|^2$
and
$|s|^2\leq2|s+t|^2+2|t|^2$,
together with $|w\cdot d|^2\leq r|d|^2$, we obtain the unifomr equivalence of the two dual metrics
\[
 \frac{1}{1+2r}|(b,d)|_{\widetilde g}^2
 \leq
 |(b,d)|_{\rm FS}^2
 \leq
 (1+r)(1+2r)|(b,d)|_{\widetilde g}^2.
\]
\end{proof}

\begin{proof}[Proof of Theorem \ref{thm-sharp}]
Set
\[
 \eta:=dw_1\wedge\cdots\wedge dw_r.
\]
We first show that the map in \eqref{eq-bergman-map} is
well defined.
For any $h\in A^2(\mathbb D^r)$, consider the
$\bar\partial$-closed form $h(w)\omega_{n-r}(z)\wedge\eta.$
By Observation \ref{lem-metric-comparison}, it suffices to
verify the $L^2$ integrability with respect to
$\widetilde g$.
Note that
\[
 \bigl|\omega_{n-r}\wedge\eta\bigr|_{\widetilde g}^2
 =
 (n-r)(1+|z|^2)^r,
 \qquad
 dV_{\widetilde g}
 =
 (1+|z|^2)^{-r}\,dV_{n-r}\,dV_w,
\]
where $dV_{n-r}
 =
 \frac{\omega_{n-r}^{n-r}}{(n-r)!}$
and $dV_w$ is the volume form of
$g_{{\rm E},r}$.
Thus we obtain 
\[
\begin{aligned}
 \bigl\|h(w)\omega_{n-r}(z)\wedge\eta
 \bigr\|_{\widetilde g}^2
 &=
 (n-r)
 \left(
 \int_{\mathbb C^{n-r}}dV_{n-r}
 \right)
 \left(
 \int_{\mathbb D^r}|h(w)|^2\,dV_w
 \right)\\
 &<+\infty,
\end{aligned}
\]
where the finiteness
$\int_{\mathbb C^{n-r}}dV_{n-r}<+\infty$
 follows from the fact  that  $\omega_{n-r}$ is the restriction of the Fubini-Study metric form on 
$\mathbb P^{n-r}$.

We now prove the injectivity of the map in \eqref{eq-bergman-map}.
Suppose that the image of $h$ under
\eqref{eq-bergman-map} is zero.
Then there exists a square-integrable
$(r+1,0)$-form $U$ such that
\begin{equation}\label{eq-dbar-U}
 \bar\partial U
 =
 h(w)\omega_{n-r}(z)\wedge\eta.
\end{equation}
By Observation \ref{lem-metric-comparison}, the form $U$
is also square-integrable with respect to
$\widetilde g$.

Note that the component of $U$ containing all the holomorphic
differentials $dw_1,\ldots,dw_r$ can be written uniquely as
\[
 \beta(z,w)\wedge\eta,
\]
where $\beta$ is a $(1,0)$-form in the $z$ variables.
Thus we can write
\begin{equation}\label{U-beta-eta+U'}
    U=\beta(z,w)\wedge\eta+U',
\end{equation}
where $U'$ is the sum of the remaining holomorphic
differential components.
By the construction of $\widetilde g$, the
$z$- and $w$-cotangent directions are orthogonal with respect
to $\widetilde g$. Consequently, the induced decomposition of
exterior powers according to the $z$- and $w$-degrees is also
orthogonal.  In particular, $\beta(z,w)\wedge\eta$ is orthogonal to $U'$ w.r.t.  $\widetilde g$.  
 Hence we obtain
 \begin{equation}\label{fini-beta-eta}
     \int_{\mathbb C^{n-r}\times\mathbb D^r}
 |\beta(z,w)\wedge\eta|_{\widetilde g}^2
 \,dV_{\widetilde g}\leq ||U||^2_{\widetilde g}
 <+\infty.
 \end{equation}

On the other hand, we also have
\[
 |\eta|_{\widetilde g}^2=(1+|z|^2)^r
\]
and
\[
 dV_{\widetilde g}
 =
 (1+|z|^2)^{-r}\,dV_{n-r}\,dV_w.
\]
Again by the  orthogonality reason, we then obtain
\[
 |\beta(z,w)\wedge\eta|_{\widetilde g}^2
 =
 (1+|z|^2)^r
 |\beta(z,w)|_{g_{{\rm FS},n-r}}^2.
\]
Consequently,
\[
 |\beta(z,w)\wedge\eta|_{\widetilde g}^2
 \,dV_{\widetilde g}
 =
 |\beta(z,w)|_{g_{{\rm FS},n-r}}^2
 \,dV_{n-r}\,dV_w.
\]
By Fubini's theorem and \eqref{fini-beta-eta}, we obtain
\[
 \int_{\mathbb D^r}
 \left(
 \int_{\mathbb C^{n-r}}
 |\beta(z,w)|_{g_{{\rm FS},n-r}}^2
 \,dV_{n-r}
 \right)dV_w
 <+\infty,
\]
and in particular, for almost every $w\in\mathbb D^r$,
\[
 \beta(\,\cdot\,,w)\in
 L^2_{1,0}(\mathbb C^{n-r},\omega_{n-r}).
\]
Taking the component of \eqref{eq-dbar-U} containing all
the holomorphic differentials $dw_1,\ldots,dw_r$ and having
antiholomorphic $z$-degree one, we obtain via \eqref{U-beta-eta+U'} that
\begin{equation}\label{eq-slice-equation}
 \bar\partial_z\beta(z,w)
 =
 h(w)\omega_{n-r}(z).
\end{equation}

In conclusion,  for almost every
$w\in\mathbb D^r$ the form
$\beta(\,\cdot\,,w)$ belongs to
$L^2_{1,0}(\mathbb C^{n-r},\omega_{n-r})$
and satisfies \eqref{eq-slice-equation} on
$\mathbb C^{n-r}$.
If $h(w)\neq0$ for such a point $w\in\mathbb D^r$, then
$h(w)^{-1}\beta(\,\cdot\,,w)$ would be an
$L^2$ solution of
\[
 \bar\partial_z
 \bigl(h(w)^{-1}\beta(\,\cdot\,,w)\bigr)
 =
 \omega_{n-r},
\]
which contradicts Observation \ref{lem-affine-obstruction}.
Thus $h(w)=0$ for almost every $w\in\mathbb D^r$.
Since $h$ is holomorphic, it follows that $h$ vanishes
identically on $\mathbb D^r$. This proves the injectivity of 
 the map in \eqref{eq-bergman-map}.

Since $A^2(\mathbb D^r)$ is infinite dimensional,  $H^{r+1,1}_{(2)}
 \bigl(\mathbb C^{n-r}\times\mathbb D^r\bigr)$
is infinite dimensional.  This completes the proof of Theorem \ref{thm-sharp}.
\end{proof}

\section*{Acknowledgements}
The author would like to express his sincere gratitude to Professors Xiaoshan Li and Sheng Rao for many valuable discussions on related topics concerning $L^2$ estimates and $L^2$ cohomology.

\end{document}